# Optimisation of power transmission systems using a discrete Tabu Search method

**AM CONNOR** and **D G TILLEY**
University of Bath, UK

**Abstract**
*This paper presents a brief description of the Tabu Search method and shows how it can be applied to two different power transmission systems. Examples are presented from two transmission systems. In the first example a mechanical transmission system is considered. A four bar mechanism is synthesised in order to produce a desired output motion. The second example is a hydrostatic transmission operating under closed loop control in order to maintain a constant operating speed as the loading conditions change.*

**NOTATION**

| | |
|---|---|
| $\omega_{desired}$ | Desired motor speed |
| $\omega_{actual}$ | Actual motor speed |
| $Q_p$ | Pump flow |
| $Q_{rv}$ | Relief valve flow |
| n | Number of simulation timesteps |

## 1 INTRODUCTION

Numerical optimisation techniques can be applied as a tool in the design of engineering systems and are particularly applicable to parametric design problems. A parametric design problem is one where the general form or type of solution is known but it is necessary to determine values for the design variables to ensure that the system produces the desired response.

In this paper two different examples of parametric design are considered. In the first example a four bar mechanism is designed so that the coupler point traces a desired profile. In the second example a hydraulic circuit is designed so that a given speed response is satisfied.

These two examples have been chosen as they represent different aspects of power transmission. The first example considers non-uniform motion generation whilst the second example requires only a rotational speed output. The purpose of this paper is to demonstrate that a single optimisation method can be applied to such different problems. This leads to a coherence in the approach for the optimisation of power transmission systems.

## 2  TABU SEARCH

Tabu Search is a metaheuristic which is used to guide local search methods towards a globally optimum solution. The power of the Tabu Search algorithm is derived from the use of flexible memory cycles of differing time spans. These memory cycles are used to control, intensify and diversify the search in order to find a suitable solution. A full description of the Tabu Search method is given in (1,2) and a description of the current implementation used is given in (3). A brief description of the main elements of the method will be given here.

### 2.1  Short term memory

The most simple implementation of a Tabu search is based around the use of a hill climbing algorithm. Once the method has located a locally optimal solution, the short term memory is used to force the search out from this optima in a different direction. The short term memory constitutes a form of aggressive search that seeks to always make the best allowable move from any given position. Short term memory is implemented in the form of *tabu restrictions* which prevent the search from cycling around a given optimal or sub-optimal position.

This concept can be illustrated by considering the diagram shown in Figure 2.1. This shows a contour plot of a two dimensional function which contains one local and one global optima and the aim of the search is to find the location with the lowest value.

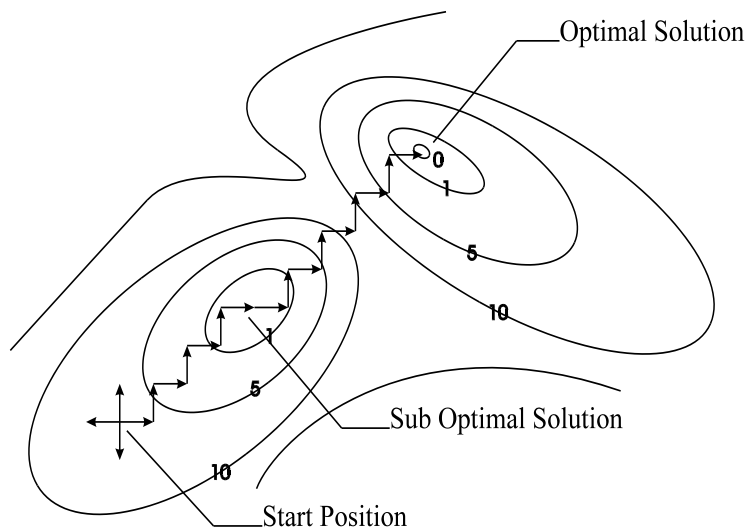

Figure 2.1 : Action of tabu restrictions

From the indicted start position the local search algorithm quickly locates the locally optimum solution without the Tabu restrictions being considered as a continuous descent is possible. However, when the search reaches the local optima the aggressive nature of the Tabu Search forces the algorithm out of the optima such that the increase in objective function is smallest. Because the last 'n' visited solutions are classed as Tabu, the search cannot leave the optima along the reverse trajectory from which it entered and once it has left the optima it cannot enter it again. The algorithm then forces the search to climb out of the local optima and locate the global optima.

## 2.2 Search intensification and diversification

The short term memory of the Tabu Search enables the method to leave locally optimum solutions in order to locate the global optimum of a function. However, short term memory alone does not ensure that the search will be both efficient and effective. Search intensification and diversification techniques are often used to first focus the search in particular areas and then expand the search to new areas of the solution space. This is normally achieved by the use of intermediate and long term memory cycles.

Intermediate and long term memory cycles generally use similar lists of previously visited solutions to guide the search. In the specific implementation used in this work the intermediate term memory cycle is based on a list of the 'm' previously visited best solutions. This list is therefore only updated when a new improved solution is found as opposed to whenever a move is made. At certain stages throughout the search process a degree of intensification is achieved by reinitialising the search at a new point generated by considering similarities between the solutions contained in the intermediate memory list.

The simplest implementation of search diversification is through a simple random refreshment of the search point though more complex methods can be implemented which utilise a long term memory cycle.

## 2.3 Hill climbing algorithm

The underlying hill climbing algorithm used in this work is based upon the method developed by Hooke and Jeeves (4). This method consists of two stages. The first stage carries out an initial exploration around a given base point. Once this exploration has been carried out and a new point found, the search is extended along the same vector by a factor $k$. This is known as a pattern move. If the pattern move locates a solution with a better objective value than the exploration point, then this point is used as a new base point and the search is repeated. Otherwise, the search is repeated using the exploration point as a new base point.

This work differs from previous work in the algorithm used to reduce the step sizes as the search progresses. Previous work reduced the step size by half when certain conditions were met. However, in this work, the step size is reduced by a fixed multiple of the minimum step size so producing an entirely discrete search.

# 3    MECHANISM SYNTHESIS

The design of mechanisms consists of two stages of synthesis. The first stage is type synthesis where the designer selects a mechanism configuration which is likely to produce the desired output motion. The second stage of the process is dimensional synthesis where the parameter values for the mechanism link lengths are selected to ensure that the output motion requirement is matched exactly.

The work presented in this paper deals only with the dimensional synthesis of four bar path generating mechanisms. Such a mechanism is shown in Figure 3.1.

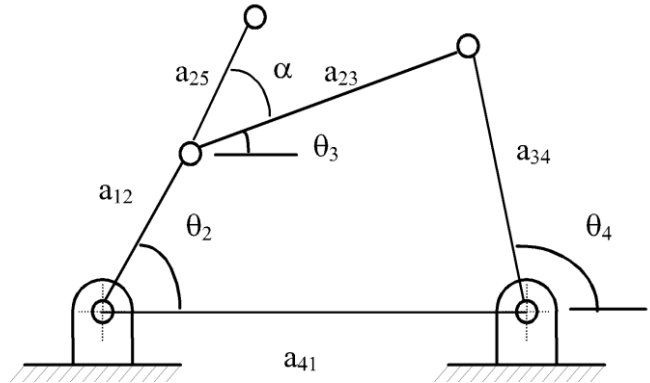

Figure 3.1 : Four bar mechanism

The mechanism consists of four links connected together by revolute joints. One link ($a_{41}$) is fixed as the ground link and in this work is constrained to be horizontal thus eliminating the angle $\theta_1$. The link $a_{12}$ is the input crank whilst link $a_{34}$ is the rocker. The coupler point is attached to the link $a_{23}$ and its position is defined by the length of the 'virtual' link $a_{25}$ and the angle $\alpha$. If the crank is able to rotate fully then this point traces a closed curve.

The condition for full crank rotation is known as the Grashof mobility criteria and can be stated as an inequality constraint;

$$l + s < p + q$$

where l and s are the lengths of the longest and shortest links respectively and p and q are the lengths of the intermediate links. This condition describes only the potential for full rotation. In order for the mechanism to be a crank-rocker type the input crank must be the shortest link and it must be adjacent to the longest link.

In this work the search is constrained to search only for Grashof crank-rocker mechanisms. The design parameters used are the lengths of the five links in the mechanism and the angle $\alpha$. Each link is constrained to certain values defined by the limits $10 \leq a_{ij} \leq 250$ and the angle $\alpha$ can take any value between -180° and 180°.

Mechanism synthesis through the use of numerical optimisation requires the definition of an objective function which describes the performance of the mechanism in relation to the desired performance. In this case a very simple objective function is used. For each step in the cycle the position of the end effector is calculated and then compared to the desired position. The error between the two points is squared so that only the magnitude of the error is incorporated into the objective function. These values are summed around the cycle of the mechanism to give the value of the objective function or a "figure of merit". One test motion is used to illustrate the effectiveness of the synthesis approach. This motion is defined by twelve precision points at which the error is calculated.

### 3.1    Test motion

The test motion used is relatively complex due to the presence of a double point in the closed curve. Previous work (5, 6, 7) which utilised a Genetic Algorithm as a synthesis method has shown it is often difficult to synthesise mechanisms for curves of this nature with a relatively low number of precision points due to the corresponding lack of definition. The results obtained from the five runs for this test motion are given in Table I.

| RUN | $a_{12}$ | $a_{23}$ | $a_{34}$ | $a_{41}$ | $a_{25}$ | $\alpha$ | OBFN | NO. EVALS |
|---|---|---|---|---|---|---|---|---|
| 1 | 30 | 52 | 57 | 71 | 50 | -52 | 1.429809 | 1725 |
| 2 | 29 | 54 | 62 | 73 | 50 | -57 | 14.88227 | 2384 |
| 3 | 30 | 53 | 51 | 70 | 50 | -45 | 13.07108 | 1511 |
| 4 | 24 | 56 | 151 | 168 | 49 | -67 | 711.3359 | 3374 |
| 5 | 30 | 53 | 55 | 70 | 50 | -50 | 2.759582 | 1540 |

Table I : Results for test motion

It can be seen that in most runs the method is locating solutions in the same region of the solution space. The large number of local optima in the solution space are preventing the method from finding the exact same solution. However, there is one solution which exhibits a large error between the actual and desired output motions (obfn ≈ 711.3). This is likely to be due to the complex motion being insufficiently defined by only twelve precision points so leading to the possibility of poor convergence.

The convergence of the best of the five runs is shown in Figure 3.2. In order to reduce the scale of the objective function values, the initial unfeasible solution has been truncated. It can be seen that the method quickly locates solutions with relatively low objective function values and converges to a near optimal solution.

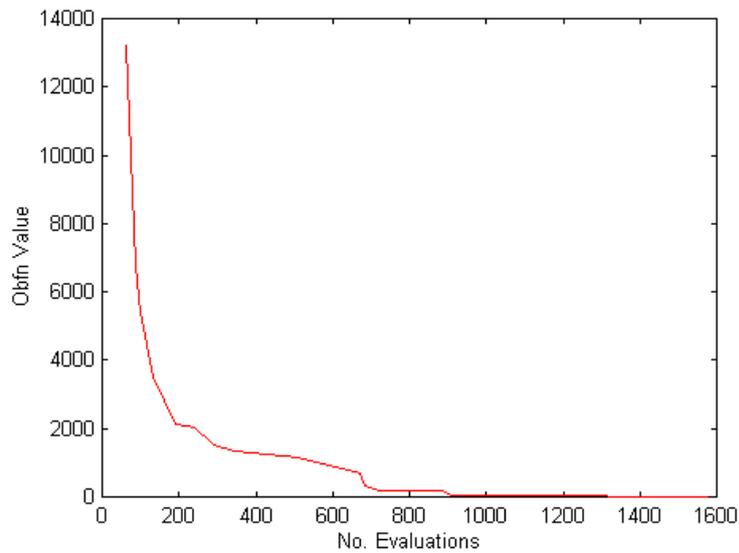

Figure 3.2 : Convergence of best solution

Figure 3.3 shows the output of the best solution in relation to the twelve precision points used to define the desired motion.

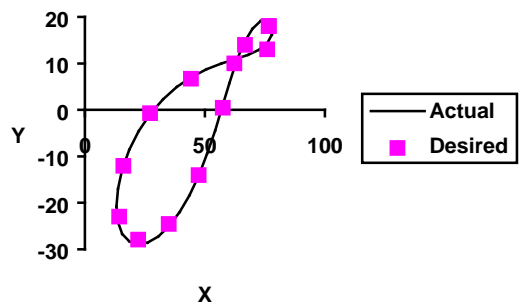

Figure 3.3 : End effector output motion

It can be seen that the method has found a solution with very low errors between the desired and actual curves.

## 4    OPTIMISATION OF HYDRAULIC SYSTEMS

In addition to applying the Tabu Search method to problems in mechanism synthesis, the Tabu Search has been applied to a number of problems in the field of fluid power system design. Figure 4.1 shows one of the test circuits used.

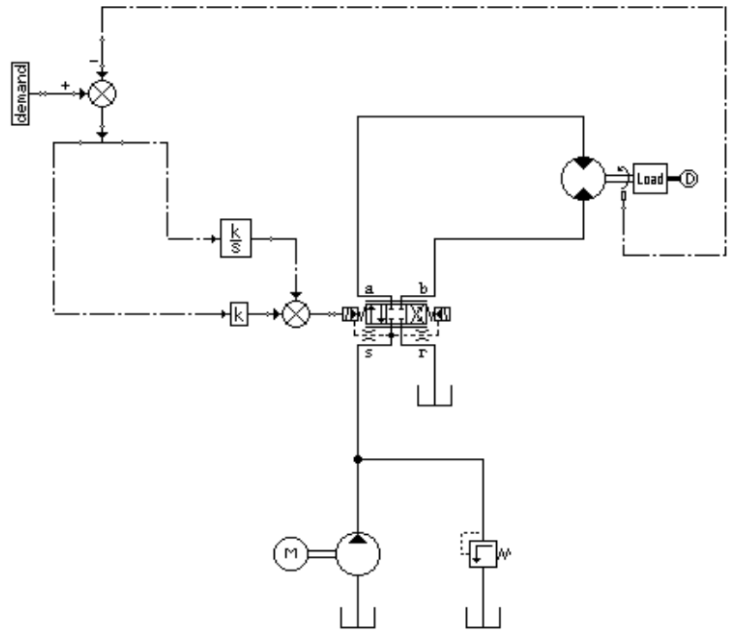

Figure 4.1 : Hydraulic test circuit

In this circuit an electrically operated proportional servo valve is use to modulate the flow from the pump so that the speed of the motor follows a desired profile. The demand for this motor is for it to ramp up to an operating speed of 300rpm in the first second of operation. The servo valve is also used to maintain a constant speed of 300 rpm despite there being a change in operating condition. After 2 seconds there is a step change in applied torque from 0 to 300Nm.

In this example the design parameters are the pump and motor displacements and the integral gain of the control system. The pump displacement is limited between the values of 10 cc/rev and 500 cc/rev. The motor displacement is limited between 50 cc/rev and 800 cc/rev whilst the integral gain can take any value between 0.01 and 50. The objective function used is described by equation 4.1.

$$obfn = \sum_{i=0}^{n}\left[(\omega_{desired} - \omega_{actual})^2 \times \left(1 + \frac{Q_{rv}}{Q_p}\right)\right]$$

… (4.1)

At each timestep of the simulation the error between the desired speed and the actual speed is calculated and this is penalised by considering the proportion of the flow produced by the pump that is returned to tank through the relief valve. The sum of these values for the whole simulation provides the objective function value.

Table II shows the results for five trial runs on this problem. It can be seen that the Tabu Search is locating solutions in the same region of the solution space for each run, though once again each solution is marginally different which suggests that the objective function is again relatively complex. However, it is important to point out that for each solution the integral gain value is equal to the upper bound.

| RUN | PUMP SIZE | MOTOR SIZE | INTEGRAL GAIN | OBFN | NO. EVALS |
|---|---|---|---|---|---|
| 1 | 150 | 740 | 50 | 1154.75 | 811 |
| 2 | 140 | 695 | 50 | 1022.14 | 973 |
| 3 | 150 | 740 | 50 | 1154.75 | 876 |
| 4 | 155 | 765 | 50 | 1231.67 | 973 |
| 5 | 155 | 765 | 50 | 1231.67 | 896 |

Table II : Results for hydraulic circuit

The convergence of the best solution obtained is shown in Figure 4.2. The method locates solutions with low objective function values and then converges to a near optimal solution.

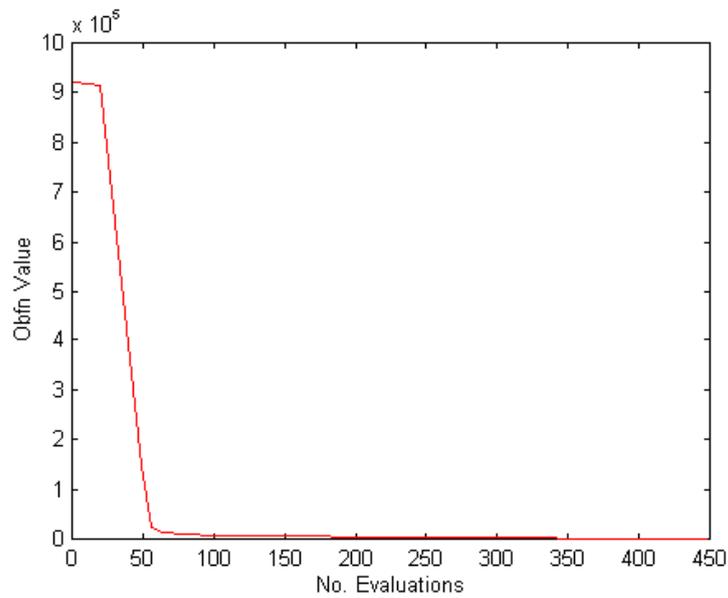

Figure 4.2 : Convergence of best solution

This solution has the speed response shown in Figure 4.3. The solution has achieved a reasonable response for the initial ramping of the speed to the desired output speed. There is an initial error that is quite large which is likely to be due to the fact that the circuit starts with the servo valve fully closed and some lag occurs as it opens to the demand position.

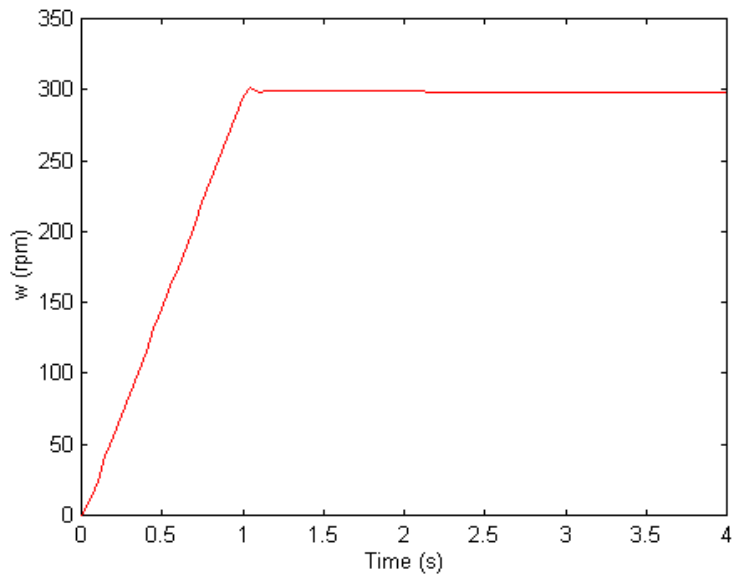

Figure 4.3 : Speed response of best solution

Some oscillation occurs when the motor reaches the desired steady state operating speed and again some disturbance can be seen to occur when the load torque is applied after two seconds of operation. However, considering the simple nature of the control system such small errors are acceptable.

However, by examining the relief valve flow throughout the first four seconds of operation it can be seen that the circuit does not exhibit entirely satisfactory performance. The relief valve flow is shown in Figure 4.4.

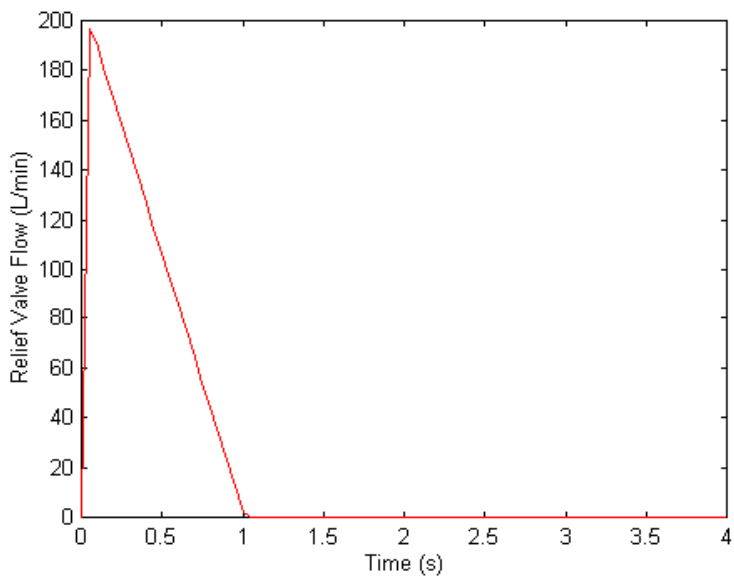

Figure 4.4 : Relief Valve Flow

Whilst the circuit has zero steady state relief valve flow, it can be seen that the valve does open in the first second of operation allowing flow to return to tank. This opening of the relief valve is caused by the required restriction in flow rate to ensure that the motor speed follows the desired ramp up to the final operating speed.

## 5  DISCUSSION

The results presented show that the Tabu Search method can be successfully applied to problems from widely differing domains. For the mechanism synthesis problem the method compares favourably to other approaches and is producing a solution that is feasible and provides a very close approximation to the desired output motion. This occurs despite the desired motion being quite complex, despite there only being twelve precision points defining the motion. This compares to the previously used Genetic Algorithm (5) which failed to located acceptable solutions for complex output motions.

For the hydraulic circuit optimisation the method has located high quality solutions with a very small number of objective function evaluations. This small number of evaluations is essential as the dynamic simulation of hydraulic circuits can be exceptionally time consuming. In both cases the Tabu Search method generally locates solutions in the same region of the solution space which differ by only small amounts. This consistency has not been apparent in previous work (3), probably due to less rigorous objective function definitions. The small variations in parameter values testify to the large number of local optima present in real optimisation problems. In future work it may be possible to eliminate these small variations by changing the implementation of the Tabu Search method memory cycles allowing it to move more easily between local optima.

## 6  CONCLUSIONS

This paper illustrates the effectiveness of the Tabu Search method on problems from two different domains. The first examples of mechanism synthesis are constrained by the mobility requirements for full crank rotatability whist in the second example of hydraulic circuit optimisation there are no constraints other than the boundaries for the design variables. In all the test cases the method has rapidly located feasible solutions which meet the desired performance characteristics so indicating that the Tabu Search is possibly a generic optimisation method which requires little or no modification in order to perform well on different problems.

## 7  ACKNOWLEDGEMENTS

The research in this paper is funded by the Engineering and Physical Sciences Research Council under grant no. 86796. This support is gratefully acknowledged.